\documentclass[12pt]{article}
\usepackage[utf8]{inputenc}
\usepackage{color}
\usepackage{amsmath,amssymb,amsthm,amsfonts}
\usepackage{float}
\usepackage{graphicx}
\RequirePackage{hyperref}

\newtheorem{remark}{Remark}

\newtheorem{definition}{Definition}
\newtheorem{theorem}{Theorem}

\usepackage{geometry}
\begin{document}
\begin{center}
\textsc{{ \Large Necessary conditions to a fractional variational problem }\\[0pt]}
\vspace{0.4cm} \hspace{0.2cm} Melani Barrios$^{1,2}$ \hspace{0.5cm} Gabriela Reyero$^{1}$ \hspace{0.5cm} Mabel Tidball$^{3}$ \\[0pt]
\end{center}

\scriptsize                                                                               
 $\,^{1}$  Departamento de Matem\'atica, Facultad de Ciencias Exactas, Ingeniería y Agrimensura, Universidad Nacional de Rosario, Avda. Pellegrini $250$, S$2000$BTP Rosario, Argentina.
 
$\,^{2}$ CONICET, Departamento de Matem\'atica, Facultad de Ciencias Exactas, Ingeniería y Agrimensura, Universidad Nacional de Rosario, Avda. Pellegrini $250$, S$2000$BTP Rosario, Argentina.

$\,^{3}$ CEE-M, Universidad de Montpellier, CNRS, INRA, SupAgro, Montpellier, France.\\

\normalsize

Correspondence should be addressed to melani@fceia.unr.edu.ar\\

\begin{abstract} 
In order to solve fractional variational problems, there exist two theorems of necessary conditions: an Euler-Lagrange equation which involves Caputo and Riemann-Liouville fractional derivatives, and other Euler-Lagrange equation that involves only Caputo derivatives. In this article, we make a comparison solving a particular fractional variational problem with both methods to obtain some conclusions about which method gives the optimal solution.
\end{abstract}




\section{Introduction}
The fractional variational calculus is a recent field, started in 1997, where classical variational problems are considered, but in the presence of fractional derivatives or integrals \cite{Agr, AlTo, Mal4}.  

In the last years numerous works have been developed tending to extend the theory of the variational calculus in order to be able to be applied to problems of fractional variational calculus.
This is fundamentally due, on the one hand, to an important development of the fractional calculus both from the mathematical point of view and its applications in other areas (electricity, magnetism, mechanics, dynamics of fluids, medicine, etc, \cite{Alm, Bas0, Goo, GoRe, Hil, Kil}), which has led to great growth in its study in recent decades.
On the other hand, the fractional differential equations establish models far superior to those that use differential equations with integer derivatives because  they incorporate into the model issues of memory \cite{FeSa} or later effects that are neglected in the models with classical derivative. 

There are several definitions of fractional derivatives  \cite{Die,Kil}. The most commonly used are the Riemann-Liouville fractional derivative and the Caputo fractional derivative. It is important to remark that while the Riemann-Liouville fractional derivatives \cite{Old} are historically the most studied approach to fractional calculus, the Caputo \cite{Cap1,Cap2} approach to fractional derivatives is the most popular among physicists and scientists, because the differential equations defined in terms of Caputo derivatives require regular initial and boundary conditions. Furthermore, differential equations with Riemann-Liouville derivatives require nonstandard fractional initial and boundary conditions that lead, in general, to singular solutions, thus limiting their application in physics and science \cite{Her,Hil}.

In order to solve fractional variational problems, there exist two theorems of optimality conditions: an Euler-Lagrange equation which involves Caputo and Riemann-Liouville fractional derivatives \cite{Agr2, Alm3, Alm,  Alm2, AlTo, BaTru, Bas,  Mal, Mal3, Mal4, OdMaTo}, and other Euler-Lagrange equation that involves only Caputo derivatives \cite{AnFe, BaRe2, BaReTid, BlaCi, Ka, LaTo}.

 
 
In the present work we will make a comparison between the solutions of the two Euler-Lagrange equations.

The paper is organized as follows: some basic definitions of fractional derivatives and fractional variational problems are shown in section two. Section three presents a particular fractional variational problem, solutions of it using both methods and the comparison between them. We end this paper with our conclusions.

\section{Mathematical tools}
\subsection{Introduction to fractional calculus}
In this section, we present some definitions and properties of the Caputo and Riemann-Liouville fractional calculus. For more details on the subject and applications, we refer the reader to \cite{Die, Old, Pod}.
\begin{definition}
The Mittag Leffler function with parameters $\alpha , \, \beta$, is defined by 
\begin{equation}
 E_{\alpha,\beta} (z) = \displaystyle \sum_{k=0}^{\infty} \frac{z^k}{\Gamma(\alpha k + \beta)}
\label{mittag}
\end{equation}
for all $z\in \mathbb{C}$.
\end{definition}
\begin{definition}
The Gamma function, $\Gamma: (0, \infty)\rightarrow \mathbb{R}$, is defined by
\begin{equation}
 \Gamma(x) = \int_{0}^{\infty} s^{x-1} e^{-s} \, ds.
\label{gamma}
\end{equation}
\end{definition}
\begin{definition}
The Riemann-Liouville fractional integral operator of order ${\alpha \in \mathbb{R}^{+}_{0}}$ is defined in $L^1[a,b]$ by
\begin{equation}
 \,_{a}I_{x}^{\alpha} [f] (x) = \dfrac{1}{\Gamma(\alpha)} \int_{a}^{x} (x-s)^{\alpha -1} f(s) \, ds.
\label{frac1}
\end{equation}
\end{definition}
\begin{definition} 
If $f \in L^1[a,b]$, the left and right Riemann-Liouville fractional derivatives of order $\alpha \in \mathbb{R}^{+}_{0}$ are defined, respectively, by
\[ \,^{RL}_{a}D_{x}^{\alpha}[f](x)= \dfrac{1}{\Gamma(n-\alpha)}\dfrac{d^{n}}{dx^{n}}\int_a^{x}(x-s)^{n-1-\alpha}f(s)ds\]
and
\[ \,^{RL}_{x}D_{b}^{\alpha}[f](x)= \dfrac{(-1)^{n}}{\Gamma(n-\alpha)}\dfrac{d^{n}}{dx^{n}}\int_x^{b}(s-x)^{n-1-\alpha}f(s)ds,\]
with $n=\left\lceil \alpha \right\rceil$.
\label{defRL} 
\end{definition}
\begin{definition}
If  $\tfrac{d^{n}f}{dx^{n}} \in L^1[a,b]$, the left and right Caputo fractional derivatives of order $\alpha \in \mathbb{R}^{+}_{0}$ are defined, respectively, by
\[ \,_{a}^{C}D_{x}^{\alpha}[f](x)= \dfrac{1}{\Gamma(n-\alpha)}\int_a^{x}(x-s)^{n-1-\alpha}\dfrac{d^{n}}{ds^{n}}f(s)ds\]
and
\[ \,_{x}^{C}D_{b}^{\alpha}[f](x)= \dfrac{(-1)^{n}}{\Gamma(n-\alpha)}\int_x^{b}(s-x)^{n-1-\alpha}\dfrac{d^{n}}{ds^{n}}f(s)ds,\]
with $n=\left\lceil \alpha \right\rceil$.
\end{definition}
Now some different properties of the Riemann-Liouville and Caputo derivatives will be seen.
\begin{remark}(Relation between the Riemann-Liouville and the Caputo fractional derivatives)\\
Considering $0<\alpha<1$ and assuming that $f$ is such that $\,^{RL}_{a} D_{x}^{\alpha}[f], \, \,^{RL}_{x} D_{b}^{\alpha}[f], \, \,_{a}^{C}D_{x}^{\alpha}[f]$ and $\,_{x}^{C}D_{b}^{\alpha}[f]$ exist, then
\[\,_{a}^{C}D_{x}^{\alpha}[f](x)=\,^{RL}_{a} D_{x}^{\alpha}[f](x)-\dfrac{f(a)}{1-\alpha} (x-a)^{-\alpha}\]
and
\[\,_{x}^{C}D_{b}^{\alpha}[f](x)=\,^{RL}_{x} D_{b}^{\alpha}[f](x)-\dfrac{f(b)}{1-\alpha} (b-x)^{-\alpha}.\]
If $f(a)=0$ then
\[\,_{a}^{C}D_{x}^{\alpha}[f](x)=\,^{RL}_{a} D_{x}^{\alpha}[f](x)\]
and if $f(b)=0$ then
\[\,_{x}^{C}D_{b}^{\alpha}[f](x)=\,^{RL}_{x} D_{b}^{\alpha}[f](x).\]
\end{remark}
\begin{remark}An important difference between Riemann-Liouville derivatives and Caputo derivatives is that, being K an arbitrary constant,
\[ \,_{a}^{C}D_{x}^{\alpha} K= 0 \ ,\hspace{1cm} \,_{x}^{C}D_{b}^{\alpha}K=0,\]
however 
\[ \,^{RL}_{a}D_{x}^{\alpha}K=\dfrac{K}{\Gamma (1- \alpha)}(x-a)^{-\alpha}, \hspace{1cm} \,^{RL}_{x}D_{b}^{\alpha}K=\dfrac{K}{\Gamma (1- \alpha)}(b-x)^{-\alpha},\]
\[ \,^{RL}_{a}D_{x}^{\alpha} (x-a)^{\alpha-1}=0, \hspace{1cm} \,^{RL}_{x}D_{b}^{\alpha}(b-x)^{\alpha-1}=0.\]
In this sense, the Caputo fractional derivatives are similar to the classical derivatives.  
\label{remark}
\end{remark}
\begin{theorem} (Integration by parts. See \cite{Kil})\\
Let $0<\alpha<1$. Let $f \in C^{1}([a,b])$ and $g \in L^1([a,b])$. Then,
\[\int_a^b g(x) \,_a^CD_x^{\alpha}f(x)\,dx=\int_a^b f(x) \,^{RL}_xD_b^{\alpha}g(x)\,dx+\left[ \,_xI_b^{1-\alpha}g(x)f(x)\right]\left|_a^b\right.\]
and
\[\int_a^b g(x) \,_x^CD_b^{\alpha}f(x)\,dx=\int_a^b f(x) \,^{RL}_aD_x^{\alpha}g(x)\,dx- \left[\,_aI_x^{1-\alpha}g(x)f(x)\right]\left|_a^b\right..\]
Moreover, if $f(a)=f(b)=0$, we have that
\[\int_a^b g(x) \,_a^CD_x^{\alpha}f(x)\,dx=\int_a^b f(x) \,^{RL}_xD_b^{\alpha}g(x)\,dx\]
and
\[\int_a^b g(x) \,_x^CD_b^{\alpha}f(x)\,dx=\int_a^b f(x)\,^{RL}_aD_x^{\alpha}g(x)\,dx.\]
\label{parts}
\end{theorem}
\subsection{Fractional variational problems}
Consider the following problem of the fractional calculus of variations which consists in finding a function $y \in \,_{a}^{\alpha}E$ that optimizes (minimizes or maximizes) the functional
\begin{equation}
J(y)= \int^{b}_{a} L(x,y,\,_{a}^{C}D_{x}^{\alpha}y) \, dx
    \label{prob}
\end{equation}
with a Lagrangian $L \in C^1([a,b]\times \mathbb{R}^2)$ and 
\[\,_{a}^{\alpha}E=\{ y: [a,b] \rightarrow \mathbb{R}: y \in C^1([a,b]), \, \,_{a}^{C}D_{x}^{\alpha}y \in C([a,b])\},\]
subject to the boundary conditions: $y(a)=y_a \, ,\, \, y(b)=y_{b}$.

Now Euler-Lagrange equations for this problem will be stated. In the first one appears both Caputo and Riemann-Liouville derivatives (Theorem \ref{teoELcrl}), meanwhile the second one only depends on Caputo derivatives (Theorem \ref{teoELcc}).

The proof of the following theorem is in \cite{Mal4}.

\begin{theorem}
If $y$ is a local optimizer to the above problem, then $y$ satisfies the next Euler-Lagrange equation: 
\begin{equation}
\dfrac{\partial{L}}{\partial y}+\,^{RL}_{x}D_{b}^{\alpha}\dfrac{\partial{L}}{\partial \,_{a}^{C}D_{x}^{\alpha}y}=0.
\label{eccrl}
\end{equation}
\label{teoELcrl}
\end{theorem}

\begin{remark} Equation (\ref{eccrl}) is said to involve Caputo and Riemann-Liouville derivatives. 
This is a consequence of the Lagrange method to optimize functionals: the application of integration by parts (Theorem \ref{parts}) for Caputo derivatives in the Gateaux derivative of the functional relates Caputo with Riemann-Liouville derivatives.
\end{remark}

\begin{remark} Equation (\ref{eccrl}) is only a necessary condition to existence of the solution. We are now interested in finding sufficient conditions. Typically, some conditions of convexity over the Lagrangian are needed. 
\end{remark}

\begin{definition} We say that $f(\underline{x},y,u)$ is convex in $S\subseteq \mathbb{R}^{3}$ if $f_{y}$ and $f_{u}$ exist and are continuous, and the condition
\[
f(x,y+y_{1},u+u_{1})-f(x,y,u) \geq f_{y}(x,y,u)y_{1}+f_{u}(x,y,u)u_{1},
\]
holds for every $(x,y,u),(x,y+y_{1},u+u_{1}) \in S.$
\end{definition}

The following theorem is valid only for the solution of the Euler-Lagrange equation involving Riemann-Liouville and Caputo derivatives (\ref{eccrl}). Its proof can be seen at \cite{AlTo}. 

\begin{theorem}
Suppose that the function $L(\underline{x},y,u)$ is convex in $[a,b]\times\mathbb{R}^{2}$.
Then each solution $y$ of the fractional Euler–Lagrange equation (\ref{eccrl}) minimizes (\ref{prob}), when restricted to the boundary conditions $y(a)=y_a$ and $y(b)=y_b$.
\label{teosuficienteRL}
\end{theorem}

 Following \cite{LaTo}, in the below theorem, we will see an Euler-Lagrange fractional differential equation only depending on Caputo derivatives. 

\begin{theorem}
Let $y$ be an optimizer of (\ref{prob}) with $L\in C^{2}\left([a,b]\times \mathbb{R}^{2}\right)$ subject to boundary conditions $y(a)=y_a \, ,\, \, y(b)=y_{b}$, then $y$ satisfies the fractional Euler-Lagrange differential equation
\begin{equation}
\dfrac{\partial{L}}{\partial y}+\,_{x}^{C}D_{b}^{\alpha}\dfrac{\partial{L}}{\partial \,_{a}^{C}D_{x}^{\alpha} y}=0.
\label{EulerLagrange}
\end{equation}
\label{teoELcc}
\end{theorem}
\begin{remark} 
We can see that the equation (\ref{EulerLagrange}) depends only on the Caputo derivatives. It is worth noting the importance that $L\in C^{2}\left([a,b]\times\mathbb{R}^{2}\right)$, without this the result would not be valid. As we remarked before, the advantage of this new formulation is that Caputo derivatives are more appropriate for modeling problems than the Riemann-Liouville derivatives and makes the calculations easier to solve because, in some cases, its behavior is similar to the behavior of classical derivatives.

From now on, when we work with the Euler-Lagrange equation that uses derivatives of Caputo and Riemann-Liouville (\ref{eccrl}), we will abbreviate it with C-RL and when we use the Euler-Lagrange equation that uses only derivatives of Caputo (\ref{EulerLagrange}), we will abbreviate it with C-C.

\end{remark}
\begin{remark} 
Unlike the equation (\ref{eccrl}), at the moment, there are not sufficient conditions for the equation (\ref{EulerLagrange}) which only involves Caputo derivatives. 
\end{remark} 

Now we present an example that we are going to solve using these two different methods, in order to make comparisons.

\section{Example}
The scope of this section is to present two different candidates to be a solution for a particular problem that arise from solving the two Euler-Lagrange equations presented in the previous section.

First, we are going to solve the classical case, where only appears an integer derivative, and then we are  going to deal with the fractional case. 

\subsection{Classical case}
The classical problem consist in finding a function $y \in \,_{a}E'$ that optimizes (minimizes or maximizes) the functional 

\[J(y)=\int_{0}^{1}\left(\left(y'(x)\right)^2-24\, y(x)\right) dx,\]
\[y(0)=0 \,,\, y(1)=0,\]

where $\,_{a}E'=\{ y: [a,b] \rightarrow \mathbb{R}: y \in C^1([a,b]) \}$.

To solve this (refer to \cite{Van}), we consider the Lagrangian 
\begin{equation}
  L(x,y,y')=\left(y'\right)^2-24\, y.
  \label{lagrangianoclasico}
\end{equation}
Its Euler-Lagrange equation is
\[\frac{\partial L}{\partial y}- \frac{\partial}{\partial x} \left(\frac{\partial L}{\partial y'}\right)=0, \]
that is, 
\[y''(x)=-12.\]
Solving this equation and taking into account that $y(0) =y(1)=0$, we obtain the solution
\begin{equation}
    y(x)=-6x^2+6x.
    \label{solclasica}
\end{equation}

\subsection{Fractional case}
The fractional problem consist in finding a function $y \in \,_{a}^{\alpha}E$ that optimizes (minimizes or maximizes) the functional 

\[J(y)=\int_{0}^{1}\left( \left(\,^{C}_{0} D^{\alpha}_{x}\left[y\right](x)\right)^2-24\, y(x)\right) dx,\]
\[y(0)=0 \,,\, y(1)=0,\]

where $\,_{a}^{\alpha}E=\{ y: [a,b] \rightarrow \mathbb{R}: y \in C^1([a,b]), \, \,_{a}^{C}D_{x}^{\alpha}y \in C([a,b]) \}$.

To solve this we consider the Lagrangian 
\begin{equation}
  L(x,y, \,^{C}_{0} D^{\alpha}_{x}\left[y\right])= \,^{C}_{0} D^{\alpha}_{x}\left[y\right]^2-24\, y.
  \label{lagrangiano}
\end{equation}

Like we said before, we are going to solve it using two methods, one with the C-RL Euler-Lagrange (\ref{eccrl}) and the other one with the C-C Euler-Lagrange equation (\ref{EulerLagrange}). 

\subsubsection{Resolution by C-RL equation}
Applying the equation (\ref{eccrl}), we obtain 
\begin{equation*}
\begin{array}{r l}
\vspace{.25cm}
\dfrac{\partial L}{\partial y} + \,^{RL}_{x} D^{\alpha}_{1}\left(\dfrac{\partial L}{\partial \,^{C}_{0} D^{\alpha}_{x}\left[y\right] }\right)&=0\\
\vspace{.25cm}
-24+\,^{RL}_{x} D^{\alpha}_{1}\left(2 \,^{C}_{0} D^{\alpha}_{x}\left[y\right]\right)&=0\\
\,^{RL}_{x} D^{\alpha}_{1}\left( \,^{C}_{0} D^{\alpha}_{x}\left[y\right]\right)&=12.
\end{array}
\end{equation*}

By definition,
\[\,^{RL}_{x} D^{\alpha}_{1}\left[(1-x)^{\beta}\right]=\frac{\Gamma (1+\beta)}{\Gamma(1+\beta-\alpha)}(1-x)^{\beta-\alpha}\]

and the property
\[\,^{RL}_{x} D^{\alpha}_{1}\left[ (1-x)^{\alpha-1} \right]=0,\]

which we have seen on remark \ref{remark}, considering $\beta=\alpha$ we can conclude 
\begin{equation}
\begin{array}{r l}
 \,^{C}_{0} D^{\alpha}_{x}\left[y\right](x)&= \dfrac{12}{\Gamma(1+\alpha)} (1-x)^{\alpha} + c_1\, (1-x)^{\alpha-1}
\end{array}
\label{change1rl}
\end{equation}

where $c_1 \in \mathbb{R}$.

Taking into account the following equalities
\begin{equation*}
    \begin{aligned}
    \vspace{.25cm}
 (-1)^n \prod_{j=0}^{n-1} (\alpha-j)&=\dfrac{\Gamma(n-\alpha)}{\Gamma(-\alpha)},\\
 (-1)^n \prod_{j=0}^{n-1} (\alpha-1-j)&=\dfrac{\Gamma(n-\alpha+1)}{\Gamma(-\alpha+1)},
    \end{aligned}
\end{equation*}

we can write
\[\begin{array}{r l}
(1-x)^\alpha&=\sum \limits_{n=0}^{\infty} \dfrac{(-1)^n \prod_{j=0}^{n-1} (\alpha-j)}{n!} x^n\\
            &=\sum \limits_{n=0}^{\infty} \dfrac{\Gamma(n-\alpha)}{\Gamma(-\alpha)} \dfrac{x^n}{n!},\\
\,\\
(1-x)^{\alpha-1}&=\sum \limits_{n=0}^{\infty} \dfrac{(-1)^n \prod_{j=0}^{n-1} (\alpha-1-j)}{n!} x^n\\
                &=\sum \limits_{n=0}^{\infty}  \dfrac{\Gamma(n-\alpha+1)}{\Gamma(-\alpha+1)} \dfrac{x^n}{n!},\\
\end{array}\]

replacing this in (\ref{change1rl}),
\[ \,^{C}_{0} D^{\alpha}_{x}\left[y\right](x)= \frac{12}{\Gamma(1+\alpha)} \sum \limits_{n=0}^{\infty} \frac{\Gamma(n-\alpha)}{\Gamma(-\alpha)} \frac{x^n}{n!} + c_1\, \sum \limits_{n=0}^{\infty}  \frac{\Gamma(n-\alpha+1)}{\Gamma(-\alpha+1)} \frac{x^n}{n!}. 
\]

Considering $\,^{C}_{0} D^{\alpha}_{x}\left[ x^{\beta}\right]= \tfrac{\Gamma (1+\beta)}{\Gamma(1+\beta-\alpha)}x^{\beta-\alpha}$ and the linearity of the Caputo derivative, we obtain

\begin{equation*}
\begin{array}{r l}
y(x)&=\dfrac{12}{\Gamma(1+\alpha)^2} x^{\alpha}  \sum \limits_{n=0}^{\infty} \dfrac{\Gamma(n+1)\Gamma(n-\alpha)\Gamma(1+\alpha)}{\Gamma(1)\Gamma(-\alpha)\Gamma(1+n+\alpha)}\dfrac{x^n}{n!}+\\
\,\\
&\qquad +\dfrac{c_1}{\Gamma(1+\alpha)} x^{\alpha} \sum \limits_{n=0}^{\infty} \dfrac{\Gamma(n+1)\Gamma(n-\alpha+1)\Gamma(1+\alpha)}{\Gamma(1)\Gamma(1-\alpha)\Gamma(1+n+\alpha)}\dfrac{x^n}{n!} +c_2,
\end{array}
\end{equation*}

where $c_2 \in \mathbb{R}$.

Using the definition of the Hypergeometric function of parameters $a$, $b$, $c$ \cite{Erl}:
\begin{equation}
\,_{2} F _{1}(a,b,c,x)=\sum \limits_{n=0}^{\infty} \frac{\Gamma(a+n)\Gamma(b+n)\Gamma(c)}{\Gamma(a)\Gamma(b)\Gamma(c+n)}\frac{x^n}{n!},
\label{hiper}
\end{equation}

we can rewrite the solution as
\[
\begin{aligned}
y(x)&=\frac{12}{\Gamma(1+\alpha)^2} x^{\alpha}\, \,_{2}F_{1}(1,-\alpha,1+\alpha,x)+\\
&\qquad +\frac{c_1}{\Gamma(1+\alpha)}x^{\alpha}\,_2F_1(1,1-\alpha,1+\alpha,x)+c_2.
\end{aligned}
\]

Taking into account that $y(0)=y(1)=0$, we obtain
\begin{equation}
\begin{aligned}
y_{RL}(x)&=\frac{12}{\Gamma(1+\alpha)^2} x^{\alpha}\, \,_{2}F_{1}(1,-\alpha,1+\alpha,x)-\\
&\qquad -\frac{6}{\Gamma(1+\alpha)^2}\frac{x^{\alpha}}{\,_2F_1(1,1-\alpha,1+\alpha,1)} \,_2F_1(1,1-\alpha,1+\alpha,x).\\
\end{aligned}
\label{solC-RL}
\end{equation}

\begin{remark}
This solution is valid only for $\alpha>0.5$ since otherwise the solution tends to infinity and does not satisfy the terminal condition.
\label{soldivergencia}
\end{remark}

Finally, since $L(x,y,u)=u^2-24y$ is a convex function, 
indeed
\[
\begin{array}{r l}
L(x,y+y_1,u+u_1)-&L(x,y,u) =(u+u_1)^2-24(y-y_1)-u^2+24y=\\
                 &=u^2+2uu_1+u_{1}^{2}-24y-24y_1-u^2+24y=\\
                 &= u_1^{2}+2uu_1-24y_1\\
                 & \geq -24y_1+2 u u_1= \partial_{2}L(x,y,u,v)y_{1}+\partial_{3}L(x,y,u,v)u_{1},
\end{array}
\]

it is verified for every $(x,y,u),(x,y+y_{1},u+u_1) \in [0,1]\times \mathbb{R}^2$, applying the theorem \ref{teosuficienteRL}, $y_{RL}$ minimizes the problem for $0.5< \alpha \leq 1$. 

\begin{remark}
We can notice that Theorem \ref{teosuficienteRL} only works for functions $y$ that satisfy the C-RL Euler-Lagrange equation, but furthermore they must satisfy the boundary conditions. In the case of not satisfying the boundary conditions (as in the case of $0<\alpha<0.5$), the theorem does not work. 
\end{remark}

\begin{remark}
We can observed that the solution (\ref{solC-RL}) tends to (\ref{solclasica}) when $\alpha$ tends to 1. This means that when $\alpha=1$, we recover the solution of the classical problem.
\end{remark}

\subsubsection{Resolution by C-C equation}
As the Lagrangian (\ref{lagrangiano}) $L \in C^{2}\left([0,1]\times \mathbb{R}^{2}\right)$, we can apply the Theorem \ref{teoELcc}. Then using the equation (\ref{EulerLagrange}), we obtain
\begin{equation*}
\begin{array}{r l}
\vspace{.25cm}
\dfrac{\partial L}{\partial y} + \,^{C}_{x} D^{\alpha}_{1}\left(\dfrac{\partial L}{\partial \,^{C}_{0} D^{\alpha}_{x}\left[y\right] }\right)&=0\\
\vspace{.25cm}
-24+\,^{C}_{x} D^{\alpha}_{1}\left(2 \,^{C}_{0} D^{\alpha}_{x}\left[y\right]\right)&=0\\
\,^{C}_{x} D^{\alpha}_{1}\left( \,^{C}_{0} D^{\alpha}_{x}\left[y\right]\right)&=12.
\end{array}
\end{equation*}

By definition,
\[\,^{C}_{x} D^{\alpha}_{1}\left[(1-x)^{\beta}\right]=\frac{\Gamma (1+\beta)}{\Gamma(1+\beta-\alpha)}(1-x)^{\beta-\alpha}\]

and the property in remark \ref{remark} that, unlike the Riemann-Liouville derivative,
\[\,^{C}_{x} D^{\alpha}_{1}\left[ d_1 \right]=0,\]

for every $d_1 \in \mathbb{R}$, considering $\beta=\alpha$ we can conclude
\begin{equation}
\begin{array}{r l}
 \,^{C}_{0} D^{\alpha}_{x}\left[y\right](x)&= \dfrac{12}{\Gamma(1+\alpha)} (1-x)^{\alpha} + d_1.
\end{array}
\label{change1c}
\end{equation}

Note that in this step this equation is different from (\ref{change1rl}), and that is why we are going to obtain two different solutions.

Now we can write
\[\frac{12}{\Gamma(1+\alpha)} (1-x)^{\alpha}= \frac{12}{\Gamma(1+\alpha)}  \sum \limits_{n=0}^{\infty} \frac{(-1)^n \prod_{j=0}^{n-1} (\alpha-j)}{n!} x^n.\]

Taking into account the following equality
$$(-1)^n \prod_{j=0}^{n-1} (\alpha-j)=\frac{\Gamma(n-\alpha)}{\Gamma(-\alpha)},$$

we obtain 
\[\frac{12}{\Gamma(1+\alpha)} (1-x)^{\alpha}= \frac{12}{\Gamma(1+\alpha)}  \sum \limits_{n=0}^{\infty} \frac{\Gamma(n-\alpha)}{\Gamma(-\alpha)} \frac{x^n}{n!}. \]

Replacing this in (\ref{change1c}),
\[\begin{array}{r l}
 \,^{C}_{0} D^{\alpha}_{x}\left[y\right](x)&=\dfrac{12}{\Gamma(1+\alpha)}  \sum \limits_{n=0}^{\infty} \dfrac{\Gamma(n-\alpha)}{\Gamma(-\alpha)} \dfrac{x^n}{n!} + d_1.
\end{array}\]

Considering $\,^{C}_{0} D^{\alpha}_{x}\left[ x^{\beta}\right]=\tfrac{\Gamma (1+\beta)}{\Gamma(1+\beta-\alpha)}x^{\beta-\alpha}$ and the  linearity of the Caputo derivative, we obtain

\[y(x)=\frac{12}{\Gamma(1+\alpha)}  \sum \limits_{n=0}^{\infty}  \frac{\Gamma(n-\alpha)}{\Gamma(-\alpha)} \frac{\Gamma(1+n)}{\Gamma(1+n+\alpha)} \frac{x^{n+\alpha}}{n!}+ d_1 x^{\alpha}+d_2,\]

where $d_2 \in \mathbb{R}$.

Using the definition (\ref{hiper}) of the Hypergeometric function of parameters $a$, $b$, $c$, we can rewrite the solution as
\[y(x)=\frac{12}{\Gamma(1+\alpha)^2} x^{\alpha} \, \,_{2}F_{1}(1,-\alpha,1+\alpha,x)+d_1 x^{\alpha}+d_2.\]

Taking into account that $y(0)=y(1)=0$, we obtain
\begin{equation}
y_{C}(x)=\frac{12}{\Gamma(1+\alpha)^2} x^{\alpha}\, \,_{2}F_{1}(1,-\alpha,1+\alpha,x)-\frac{6}{\Gamma(1+\alpha)^2}{x^{\alpha}}.
\label{solC-C}
\end{equation}

\begin{remark}
Unlike $y_{RL}$ in (\ref{solC-RL}), $y_{C}$ is valid for every $0<\alpha\leq 1$. However, we can not ensure that it is a minimum of the problem because there are no sufficient conditions theorem for the C-C Euler-Lagrange equation. 
\end{remark}

\begin{remark}
We can observe that the solution (\ref{solC-C}) tends to (\ref{solclasica}) when $\alpha$ tends to 1. This means that both solutions of each Euler-Lagrange equations, $y_{RL}$ and $y_{C}$ tend to the solution of the classical Euler-Lagrange equation when $\alpha=1$.
\end{remark}

\subsection{Comparison between methods}
In this section we are going to show some graphics in order to compare the solutions obtained from the different methods. 

\begin{figure}[H]
    \centering
    \includegraphics[width=.45\textwidth]{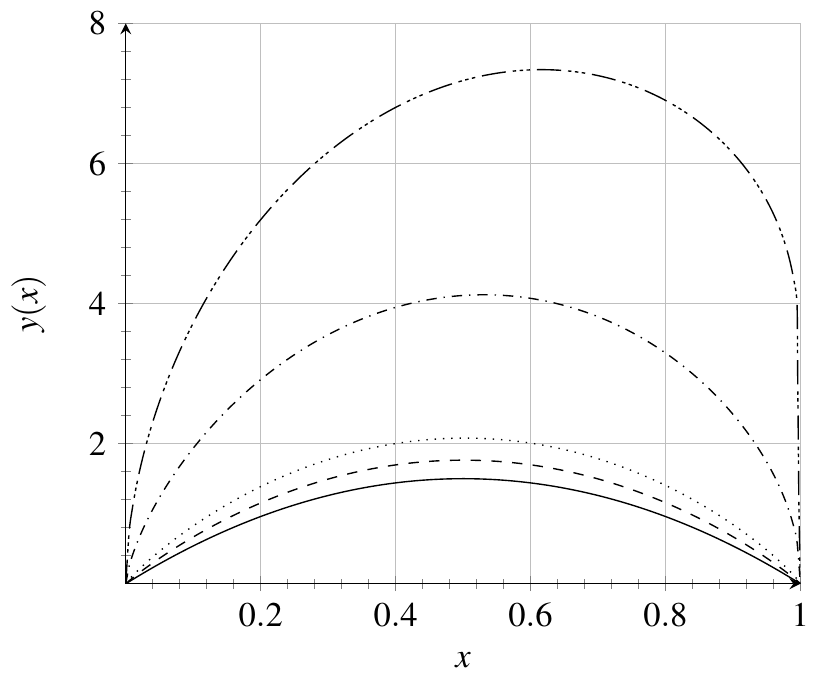} \includegraphics[width=.45\textwidth]{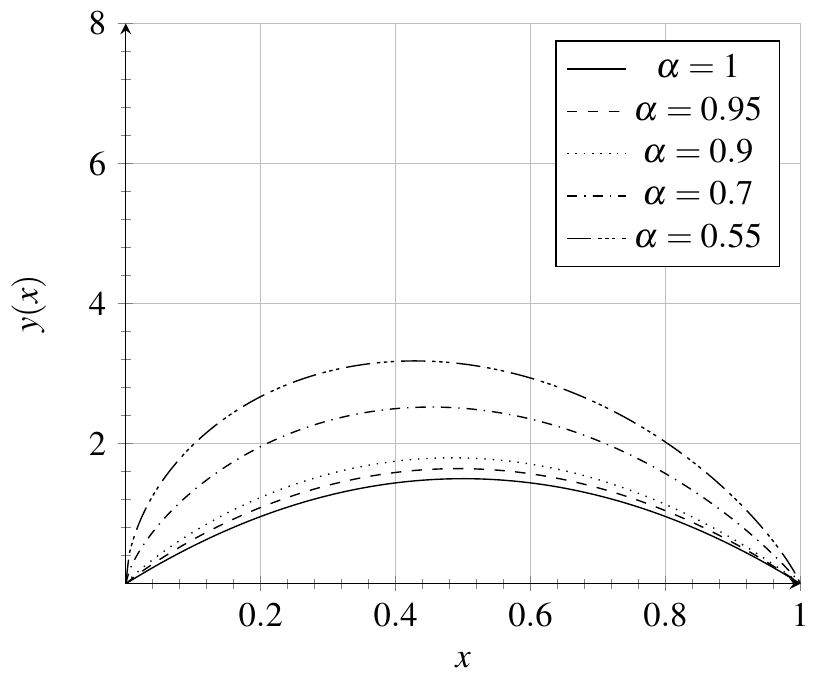}
    \caption{Convergence of $y_{RL}$ and $y_{C}$ solutions}
    \label{figconvergence}
\end{figure}

Figure \ref{figconvergence} presents the convergence of both solutions $y_{RL}$ (\ref{solC-RL}) in the left and $y_C$ (\ref{solC-C}) in the right, when we take limit as $\alpha$ approaches one. We can see that both converge to the classical solution $y$ (\ref{solclasica}). 

\begin{remark}
This figure shows us the difference between the shapes of the solutions obtained from the different methods. We can clearly see how the shapes of the solutions $y_{C}$ are more similar to the classical solution in contrast to the shapes of the solutions $y_{RL}$.
\end{remark}

\begin{figure}[t]
    \centering
    \includegraphics[width=.45\textwidth]{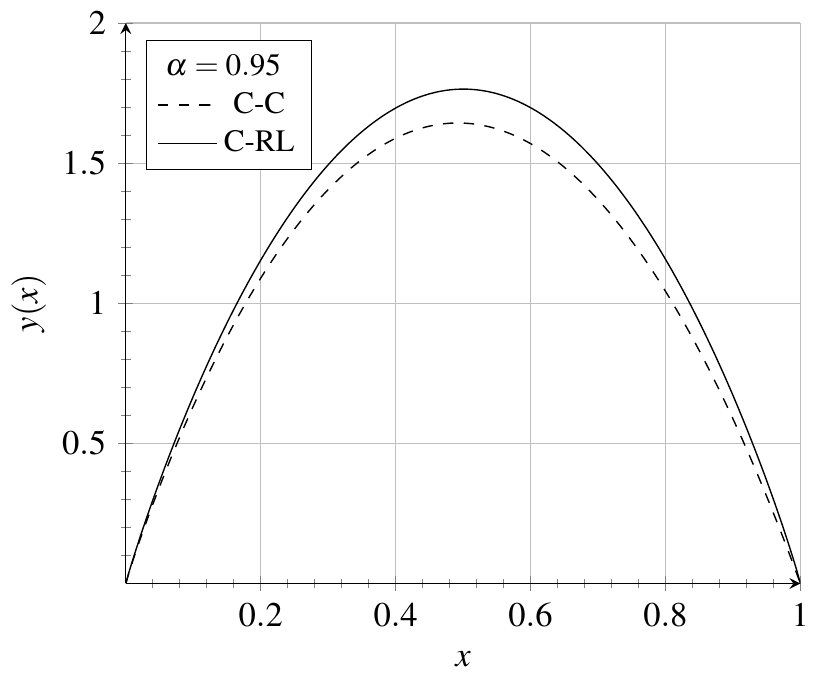} \includegraphics[width=.45\textwidth]{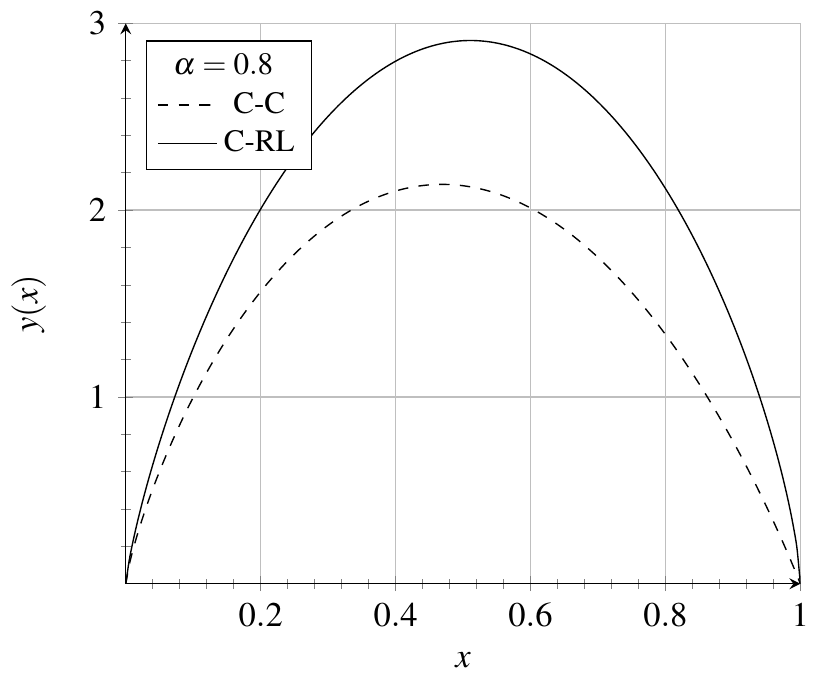}
    
    \includegraphics[width=.45\textwidth]{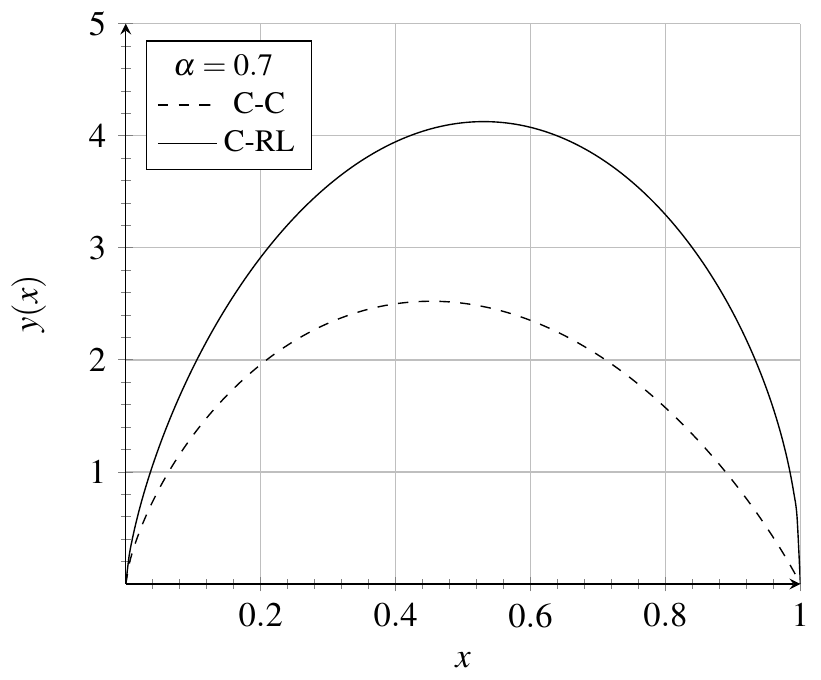} \includegraphics[width=.45\textwidth]{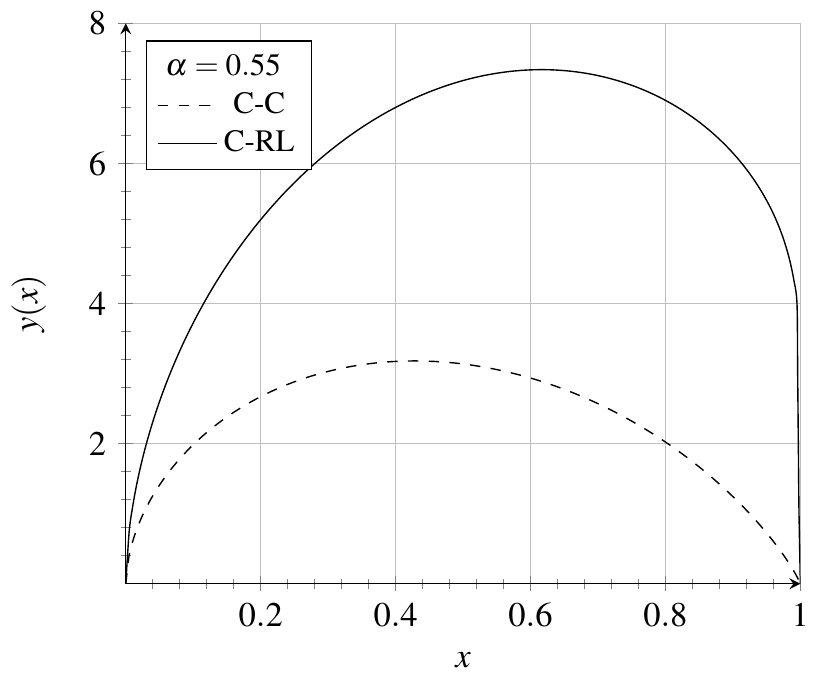}
    \caption{Comparison of the C-RL and C-C solutions}
    \label{figcomparison}
\end{figure}

Figure \ref{figcomparison} presents a comparison between both solutions $y_{RL}$ (\ref{solC-RL}) and $y_C$ (\ref{solC-C}), for different values of $\alpha$. 

\begin{remark}
In this figure we can see how the difference between the shapes of both solutions becomes more remarkable when $\alpha$ approaches 0.5, where the solutions $y_{RL}$ diverge as we saw in Remark \ref{soldivergencia}.
\end{remark}

\begin{figure}[H]
    \centering
    \includegraphics[width=.45\textwidth]{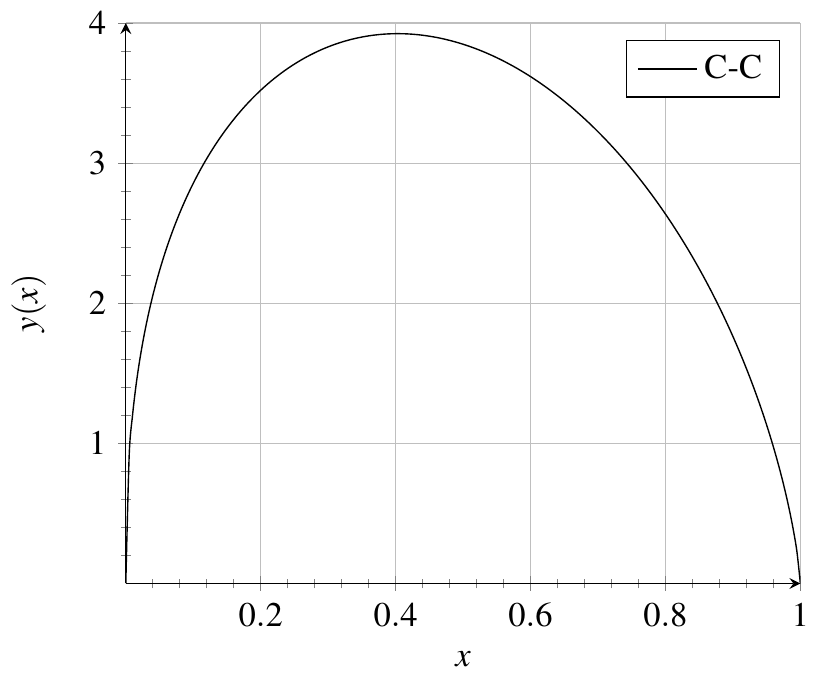}
    \caption{Solution C-C for $\alpha=0.4$}
    \label{figsolCalfa04}
\end{figure}

Figure \ref{figsolCalfa04} presents the solution $y_C$ (\ref{solC-C}) for $\alpha=0.4$. 

\begin{remark}
While the C-RL Euler-Lagrage equation does not provide us solutions for the cases $0<\alpha \leq 0.5$, the C-C equation does.
\end{remark}

Table \ref{table:minimum} presents the values obtained in each case. To calculate the integrand we approximate the Caputo fractional derivatives of both $y_{C}$ and $y_{RL}$. For this we use a method of $L_1$ type that can be seen in \cite{BaDiScTru, LiZe}. This method consists of making a regular partition of the interval $[0,1]$ as ${0=x_0 \leq x_1 \leq }$ ... $ \leq x_m=1 $, of size $h>0 $ sufficiently small, and then approximating the Caputo derivative as follows:

\[\,^{C}_{0} D^{\alpha}_{x}\left[y\right](x_{m}) =\displaystyle \sum_{k=0}^{m-1} b_{m-k-1}(y(x_{k+1})-y(x_k)),\]

where 
$$b_k=\frac{h^{-\alpha}}{\Gamma(2-\alpha)}\left[(k+1)^{1-\alpha}-k^{1-\alpha}\right].$$

Then, to calculate the integrals, we use the Riemann sums approximation.

\begin{table}[H]
\begin{center}
\begin{tabular}{| c | c | c |}
\hline
$\alpha$ & C-RL & C-C \\ \hline
1 & -12.1752 & -12.1752 \\
0.95 & -16.4431 & -14.3133 \\
0.9 & -17.3685 & -16.7006 \\
0.8 & -36.6555 & -22.2567 \\
0.7 & -60.2608 & -28.9016 \\
0.55 & -127.9983 & -40.9804 \\
0.4 & the solution does not exist & -55.5863 \\ \hline
\end{tabular}
\caption{Values obtained with C-RL and C-C \label{table:minimum}}
\end{center}
\end{table}

\begin{remark}
In Table \ref{table:minimum} we can see that as we get closer to $ \alpha=0.5$, the difference between the values is very large, being the minimum the solution of the C-RL equation, while for values $0< \alpha \leq 0.5$ obviously we only have the solution of the C-C equation, since it is the only one that verifies the border conditions.
\end{remark}

\section{Conclusions}
 In this article, two theorems of necessary conditions to solve fractional variational problems were studied: an Euler-Lagrange equation which involves Caputo and Riemann-Liouville fractional derivatives (C-RL), and other Euler-Lagrange equation that involves only Caputo derivatives (C-C).
 A particular example was presented in order to make a comparison between both conditions. 
 
 We were able to get several conclusions. The first thing is that we were able to verify that for  $0.5< \alpha \leq 1$, the minimum was obtained from the solution of the C-RL Euler-Lagrange equation, as suggested by the Theorem \ref{teosuficienteRL}. Now, for $0<\alpha \leq 0.5$, the C-RL Euler-Lagrange equation did not provide us with a solution, while C-C equation did. However, we wonder, can we ensure that the solution of the C-C equation ($y_C$) is the optimal solution for the problem at least for these values of $\alpha$? The answer to this question is NO. Although in \cite{LaTo} it was shown that the solution $y_C$ is a critical solution of the problem, and we also saw that its shape is graphically more similar to the shape of the classical solution, we cannot ensure that $y_C$ is an optimal solution or even for the cases in which $0<\alpha \leq 0.5$. If there was a Theorem of Sufficient Conditions for the C-C equations, our example would not verify it, because if it did,
this Theorem would be in contradiction with the Theorem of Sufficient Conditions for the C-RL equations (Theorem \ref{teosuficienteRL}). This means that the convexity conditions over the Lagrangian does not reach to obtain sufficient conditions for the C-C equations. Then, if there were other conditions and such a Theorem existed, in order not to contradict Theorem \ref{teosuficienteRL}, it should depends on the value of $\alpha$.

In other hand, we can observe that when $\alpha=1$, $y_C$ was the classical solution and it was the minimum of the problem, but when $\alpha$ decreased, when $0.5< \alpha<1$, these solutions were not minimums, because the $y_{RL}$ solutions were. There is a discontinuity in the $y_C$ solution when $\alpha$ goes to 1. Then, why $y_C$ would be the minimum solutions when  $0<\alpha<0.5$? If they were minimums, it would exist another discontinuity of these solutions with respect to $\alpha$.  
 
 In conclusion, while working with C-C equations make the work easier when it comes to calculations, many times we have to be careful with the implementation of this method since C-RL Euler-Lagrange equations are the ones that truly provide us with the optimal solution.




\section*{Acknowledgments}
 This work was partially supported by Universidad Nacional de Rosario through the projects ING568 ``Problemas de Control \'Optimo Fraccionario''. The first author was also supported by CONICET through a PhD fellowship.

\bibliographystyle{acm}
\bibliography{biblio} 
\end{document}